\documentclass[11pt]{article}
\usepackage{amsmath}
\usepackage{amssymb,amsmath}
\def\C{\mathop{\bf C\kern 0pt}\nolimits}
\def\DD{\mathop{\bf D\kern 0pt}\nolimits}
\def\K{\mathop{\bf K\kern 0pt}\nolimits}
\def\N{\mathop{\bf N\kern 0pt}\nolimits}
\def\Q{\mathop{\bf Q\kern 0pt}\nolimits}
\def\R{\mathop{\bf R\kern 0pt}\nolimits}
\def\T{\mathop{\bf T\kern 0pt}\nolimits}

\newcommand{\beq}{\begin{equation}}
\newcommand{\eeq}{\end{equation}}
\newcommand{\ben}{\begin{eqnarray}}
\newcommand{\een}{\end{eqnarray}}
\newcommand{\beno}{\begin{eqnarray*}}
\newcommand{\eeno}{\end{eqnarray*}}

\newtheorem{thm}{Theorem}[section]
\newtheorem{lem}{Lemma}[section]
\newtheorem{rmk}{Remark}[section]

\renewcommand{\theequation}{\thesection.\arabic{equation}}

\begin{document}
\title {Backward Uniqueness of Kolmogorov Operators}

\author{Wendong WANG$^{1,3}$ and Liqun ZHANG$^2$\\
 {$^1$ School of  Mathematical Sciences, Dalian University of Technology}\\
 { Dalian, 116024, P.R. China;}\\
  { $^2$ Institute of Mathematics, AMSS, Academia Sinica,
    Beijing, 100190}\\
  {$^3$ The Institute of Mathematical Sciences, CUHK, Hong Kong.}
\thanks{ The first author is supported by "the Fundamental Research Funds for the Central Universities" and partially by IMS of CUHK, and L. Zhang is supported by the
Chinese NSF under grant 10325104 and National Basic Research Program
of China under grant 2011CB808002. Email: wendong@dlut.edu.cn
and lqzhang@math.ac.cn}}

\date{Oct 26, 2012}
\maketitle

\begin{abstract}
The backward uniqueness of the Kolmogorov operator
$L=\sum_{i,k=1}^n\partial_{x_i}(a_{i,k}(x,t)\partial_{x_k})+\sum_{l=1}^m x_l\partial_{y_l}-\partial_t$, was proved
in this paper. We obtained a weak Carleman inequality via
Littlewood-Paley decomposition for the global backward uniqueness.
Moreover, a monotonicity inequality was also proved for the Kolmogorov equation.
\end{abstract}
{\small keywords: Carleman inequality, Kolmogorov operator,
backward uniqueness, Littlewood-Paley decomposition}

\renewcommand{\theequation}{\thesection.\arabic{equation}}
\setcounter{equation}{0}

\section{Introduction}

The Kolmogorov equation has many applications in various models (see
\cite{PP}, \cite{WZ}), for example, Prandtl's boundary layer
equation in the Crocco variables and Boltzmann-Landau equation. One
of the simplest form of the Kolmogorov operator is given in the
following equation
$$\partial_{xx} u +x\partial_y u-\partial_t
u=0.$$

In this paper, we consider the following more general backward Kolmogorov operator:
$$ Lu = (\sum_{i,k=1}^n\partial_{x_i}(a_{ik}(x,t)\partial_{x_k})+\sum_{l=1}^m x_l\partial_{y_l}-\partial_t)
u,$$ where $m\leq n$, $x=(x_1,\cdots,x_n)$,  $y=(y_1,\cdots,y_m)$, and $(x,y,t)\in{\Bbb R}^{n}\times {\Bbb R}^{m}\times (0, T)$. Assume that
the coefficients are symmetric and uniformly elliptic:
\ben \label{equ00}
a_{ik}(x,t)=a_{ki}(x,t),\quad 1\leq i,k\leq n;\quad \lambda^{-1}|\xi|^2\leq \sum_{i,k=1}^n\xi_ia_{ik}(x,t)\xi_k\leq \lambda|\xi|^2, \quad
\een
for any $(x,t)\in {\Bbb R}^{n}\times (0, T)$ and $\xi\in {\Bbb R}^{n} $. Here $\lambda>1$ is a constant.

Our interest in the backward uniqueness of Kolmogorov operator
arises from the study of regularity of Kolmogorov operators and
recent progress in the backward uniqueness of the parabolic
equations where some important applications have been found.
In fact, the backward uniqueness of parabolic operator
 in half-space is crucial for the proof of smoothness of solutions of
 Navier-Stokes equations in $L^{3,\infty}$
(see \cite{ESS2}). Some new techniques are developed in the proof of
the their result, as well as in the sequel papers (for example
\cite{EKPV}). Their main results state as follows:

Suppose that $u$ and the generalized derivatives $\partial_t u$ and
$\nabla^2 u$ are square integrable over any bounded domain of
${\Bbb R}^N_+\times (0, T)$
$$
|\Delta u+\partial_t u|\leq M (|u|+|\nabla u|),\quad |u|\leq
e^{M|x|^2}\quad {\rm in} \quad {\Bbb R}^N_+\times (0, T)$$ and $u(x,0)=0 $
in ${\Bbb R}^N_+.$ Then $u(x,t)\equiv 0,$ in $ {\Bbb R}^N_+\times (0, T).$ (see
\cite{ESS})

There is a long history of this type of backward uniqueness of
parabolic equations. In the papers like \cite{LPR} and \cite{LA},
some type of Carleman inequalities are obtained under the $C^2$
smoothness assumptions of the coefficients. On the other hand, the
well known example of Miller \cite{MI} (where an operator having
coefficients which are H\"{o}lder-continuous of order 1/6 with
respect to t and $C^{\infty}$ with respect to x does not have the
uniqueness property) shows that a certain amount of regularity
assumptions on the $a_{ik}$'s are necessary for the uniqueness.

The main idea in the proof of the backward uniqueness is to obtain
certain type of Carleman inequality, which is also useful for unique
continuation (see \cite{LIN}, \cite{CHE}, \cite{ESC}, \cite{EV},
\cite{EF}, \cite{FER}, \cite{ESS}, \cite{SP} and so on). On the other hand,
monotonicity inequality of frequency functions could also be used to
prove unique continuation , for example, see \cite{PO} and
\cite{EKPV}.

For the backward uniqueness of Kolmogorov operator, an uniform
Carleman inequality is necessary. However, it is difficult to
obtain such inequalities because of degeneracy of the Kolmogorov
operator. By combining the Littlewood-Paley decomposition and the approach  of Carleman-type inequality in Escauriaza, Seregin, and
\u{S}ver\'{a}k \cite{ESS}, we obtained a weak type of Carleman
inequality which implies the backward uniqueness property.

Our main idea is first to establish the Carleman-type inequality for
the low frequency part in the degenerated direction. Under the
assumption of coefficients which are independent of $y$, we make use of the
Littlewood-Paley decomposition of the solution $u$ in $y$ direction to
prove that the $L^2$ norm of $\triangle_j
\partial_yu$ can be controlled by the $L^2$ norm of $\triangle_j u$.
Hence the vanishing property of $\triangle_j u$ for any $j$ implies
that $u$ vanishes.

We assume that $u$ satisfies
\begin{equation} \label{equ01}
\left\{
\begin{array}{lll}Lu = (\sum_{i,k=1}^n\partial_{x_i}(a_{ik}(x,t)\partial_{x_k})+\sum_{l=1}^m x_l\partial_{y_l}-\partial_t)
u
=c(x,t)u+\sum_{i=1}^nd(x,t)_{i}\partial_{x_i} u,\\
u(x,y,0)=0\quad (x,y)\in {\Bbb R}^{n+m},\\
u, \nabla_x u, (\sum_{l=1}^m x_l\partial_{y_l}+\partial_t) u\in L^2({\Bbb R}^{n+m}\times
(0,T)),
\end{array}\right.
\end{equation}
where the coefficients are some suitable regular functions.

Let $g=\chi(t)u$, where $\chi(t)$ is a $C^{\infty}$ smooth function
for $0<t_1<t_2<T$ and
$$
\chi(t)=\left\{ \begin{array}{ll}1,\quad & \quad t\leq t_1,\\
0, \quad & \quad t> t_2. \end{array}\right.
$$
For fixed $\alpha_0>1$, we choose $j_0=\max\{j\in {\Bbb Z};$ $
4^{j+1}\leq \frac{\alpha_0-1}{4t}, t\in (0,t_2]\}$. When
$\alpha>\alpha_0$, we then choose a small $t_2$ and obtain the
following Carleman-type inequality
\begin{equation}\label{02}
\int_{{\Bbb R}^{n+m}\times (0,T)}t^{-2\alpha+1}| L \Delta_j g|^2dxdy\geq \int_{{\Bbb R}^{n+m}\times (0,T)}
\big(\frac{1}{4\lambda}t^{-2\alpha}|\partial_x \Delta_j g|^2+\frac{\alpha-1}{4}t^{-2\alpha-1}|\Delta_j g|^2\big)dxdy,
\end{equation}
where $\Delta_j$ is Littlewood-Paley decomposition operator, $j\in
{\Bbb Z}$, $j\leq j_0$. ( more details see Section 2.)

The backward uniqueness is proved by applying the above
Carleman-type inequality.

\begin{thm}\label{result01}
Suppose that $u$ satisfies ($\ref{equ01}$) and ($\ref{equ00}$). We assume that
 $a(x,t)\in C^{0,1}$ and $c(x,t)$, $d(x,t)$ are bounded functions.
Then $u\equiv 0$ in ${\Bbb R}^{n+m}\times(0,T)$.
\end{thm}

\begin{rmk}
Here $a(x,t)\in C^{0,1}$ means that $\nabla_x u$ and $\partial_t u$ are bounded. Under the condition of ($\ref{equ01}$), the operator $L$ satisfies the 
well-known H\"{o}rmander finite rank condition. 
The Kolmogorov operator, although degenerated in some sense,
still retains most of the properties of the parabolic operator. For
example, the interior regularity of the Kolmogorov operator is
similar to that of the parabolic operator (see \cite{PP}, \cite{Zhang} and \cite{WZ}), as well as the backward
uniqueness at least under some additional assumptions.
\end{rmk}

\begin{rmk}
For the general Kolmogorov operator, for example, $L_1=\partial_{xx}+x\partial_y+y\partial_z+\partial_t$, we still don't know if it has the backward uniqueness property, which is interesting and need new idea. 
\end{rmk}

We also give an alternative proof for the above results by frequency
functions as the parabolic case given by \cite{PO} and \cite{EKPV}.

\setcounter{equation}{0}
\section{ Proof of the Main Theorem}

We first introduce some of the notations which are used throughout
this paper.

Set $g=\chi(t)u$,  and $\chi(t)$ is a  $C^{\infty}$ smooth
function satisfying
$$
\chi(t)=\left\{ \begin{array}{ll}1,\quad & \quad t\leq t_1,\\
0, \quad & \quad t> t_2. \end{array}\right.
$$
where $0<t_1<t_2<T$, to be chosen.

Let $\phi(t)=(t+b)^{-\alpha}$ and $f=\phi g$, where $\alpha>0$ and
$b$ is a constant satisfying $0<b\leq t_2$.

We introduce the  Littlewood-Paley decomposition on  ${\Bbb R}^m$. Let
$\varphi(\xi)$ be a smooth cut-off function such that
$$
\varphi(\xi)=\left\{ \begin{array}{ll}1,\quad & \quad |\xi|\leq \frac12,\\
0, \quad & \quad |\xi|> 1. \end{array}\right.
$$
Let $\psi(\xi)=\varphi(\xi/2)-\varphi(\xi)$. For any integer $j$, as
usual, we denote $\triangle_j$ and $S_j$,
 $$
\triangle_j h(x)= {\mathcal F}^{-1}(\psi(\frac{\xi}{2^j}){\mathcal
F}(h)(\xi)),
 $$
 $$
S_j h(x)={\mathcal F}^{-1}(\varphi(\frac{\xi}{2^j}){\mathcal
F}(h)(\xi)),
$$ where $h(x)\in {\mathcal S}'({\Bbb R}^m)$ and ${\mathcal F}^{-1}$ is the inverse of Fourier transformation.
For $h\in L^2({\Bbb R}^m)$, denote ${\mathcal F}(h)=\hat{h}$, and we have
$$\int_{{\Bbb R}^m}|h|^2 dy=\int_{{\Bbb R}^m}|\hat{h}|^2 dy=\int_{{\Bbb R}^m}|\sum_{j\in  Z}\Delta_j h|^2 dy. $$
And it is easy to see that there exists $K>0$ only dependant on the
dimension, such that
\begin{equation}\label{equ23}
K^{-1}\int_{{\Bbb R}^m}|\sum_{j\in  Z}\Delta_j h|^2 dy  \leq \int_{{\Bbb R}^m}\sum_{j\in  Z}|\Delta_j h|^2 dy
\leq    K\int_{{\Bbb R}^m}|\sum_{j\in  Z}\Delta_j h|^2 dy\end{equation}

Let $\Omega_T={\Bbb R}^{m+n}\times (0,T)$. By our assumption, we may assume
that for some positive constant $\lambda>1$, for all $(x,t)\in {\Bbb R}^{n}\times (0,T)$ and $\xi\in {\Bbb R}^{n}$
 $$\lambda^{-1}|\xi|^2\leq \xi_ia_{ij}\xi_j\leq \lambda|\xi|^2$$ and  $$|\nabla_x a_{ij}(x,t)|, |\partial_ta_{ij} (x,t)|, |c(x,t)|, |d(x,t)|\leq
 \lambda.$$

We make the Littlewood-Paley
decomposition in the $y$ direction.
For convenience we set $$f_j=\phi \Delta_j g.$$ We are going to
prove a Carleman-type inequality for the function $\Delta_j g$ which
enable us to overcome the difficulty of the degeneracy in the $y$
direction.

\begin{lem}
Under the assumptions of Theorem \ref{result01}, the function $\Delta_j g$
satisfies the following Carleman inequality
\begin{eqnarray}\label{equ22}
\int_{\Omega_T}(t+b)^{-2\alpha+1}| L \Delta_j g|^2& \geq &
\int_{\Omega_T}\frac{1}{4\lambda}(t+b)^{-2\alpha}|\partial_x
\Delta_j g|^2+ \frac{\alpha-1}{4}(t+b)^{-2\alpha-1}|\Delta_j
g|^2.\nonumber
\end{eqnarray}
\end{lem}
\begin{rmk}
The function $\Delta_j g$ is the Littlewood-Paley decomposition in
$y$ direction of the solution $u$ of problem (\ref{equ01}). One can easily
check that
$$
\Delta_j g=\chi(t){\mathcal F}^{-1}(\psi(\frac{\xi}{2^j}){\mathcal
F}(u)(\xi)).
$$
\end{rmk}
{\it Proof:} We need to estimate the integral of the function
$\Delta_j g$ and its derivative in the $x$ direction in terms of
$$I\equiv\int_{\Omega_T}(t+b)|\phi L \Delta_j g|^2dxdydt.$$

Recall the notation $\phi(t)=(t+b)^{-\alpha}$,  $g=\chi(t)u$ and let $f_j=\phi \Delta_j g$. By
 the equation of (\ref{equ01}), on the one hand, we have
\begin{eqnarray}\label{equ24}
& &\int_{\Omega_T}(t+b)|\phi L \Delta_j g|^2dxdydt\\ &=&\int_{\Omega_T}(t+b)|(\sum_{i,k=1}^n\partial_{x_i}(a_{ik}(x,t)\partial_{x_k})
+\sum_{l=1}^mx_l\partial_{y_l}+\partial_t-\frac{\phi'}{\phi})(\phi \Delta_j g)|^2dxdydt\nonumber\\
&\geq & \int_{\Omega_T}2(t+b)(\sum_{i,k=1}^n\partial_{x_i}(a_{ik}(x,t)\partial_{x_k})-\frac{\phi'}{\phi})f_j(\sum_{l=1}^mx_l\partial_{y_l}+
\partial_t)f_j dxdydt\quad  \nonumber\\
&=&-\int_{\Omega_T}2(t+b) \sum_{i,k=1}^n a_{ik}\partial_{x_i} f_j \partial_{y_k} f_jdxdydt+\int_{\Omega_T}\sum_{i,k=1}^n \big(a_{ik}+(t+b)\partial_t a_{ik}\big)
\partial_{x_i} f_j \partial_{x_k} f_j dxdydt
\nonumber
\end{eqnarray}
where we have used the symmetry property of $a_{ik}$, and $y_k=0$ for $k>m$.

On the other hand, the usual decomposition gives another lower bound of $I$.
Let  $\tilde{\phi}=(t+b)^{\frac12}\phi=(t+b)^{\frac12 -\alpha},$
and
$\tilde{f_j}=\tilde{\phi} \Delta_j g=(t+b)^{\frac12}f_j.$ Then we
have
\begin{eqnarray}\label{equ25}
I &=&\int_{\Omega_T}|(\sum_{i,k=1}^n\partial_{x_i}(a_{ik}(x,t)\partial_{x_k})
+\sum_{l=1}^mx_l\partial_{y_l}+\partial_t-\frac{\tilde{\phi}'}{\tilde{\phi}})(\tilde{\phi} \Delta_j g)|^2dxdydt\nonumber\\
&\geq & \int_{\Omega_T}2(\sum_{i,k=1}^n\partial_{x_i}(a_{ik}(x,t)\partial_{x_k})-\frac{\tilde{\phi}'}{\tilde{\phi}})\tilde{f_j}(\sum_{l=1}^mx_l\partial_{y_l}+
\partial_t)\tilde{f_j} dxdydt\quad \nonumber \\
&=&\int_{\Omega_T}\big(-2 \sum_{i,k=1}^n a_{ik}\partial_{x_i} \tilde{f}_j \partial_{y_k} \tilde{f}_j+\sum_{i,k=1}^n\partial_t a_{ik}
\partial_{x_i} \tilde{f}_j \partial_{x_k} \tilde{f}_j +
(\frac{\tilde{\phi}'}{\tilde{\phi}})'|\tilde{f_j}|^2 \big)dxdydt\nonumber \\&
=& -\int_{\Omega_T}2(t+b) \sum_{i,k=1}^n a_{ik}\partial_{x_i} f_j \partial_{y_k} f_jdxdydt+\int_{\Omega_T}\sum_{i,k=1}^n (t+b)\partial_t a_{ik}
\partial_{x_i} f_j \partial_{x_k} f_j dxdydt
\nonumber\\
&& +\int_{\Omega_T}\frac{\alpha-\frac12}{(t+b)}|f_j|^2dxdydt
\end{eqnarray}
where in the last equality we used the fact
$\tilde{f_j}=(t+b)^{\frac12}f_j$.

Now we choose $t_2\leq (16\lambda^2)^{-1}$ and then fixed from now
on in the definition of function $\chi(t)$. Combining the two inequalities
($\ref{equ24}$) with ($\ref{equ25}$), a simple calculation yields that
\begin{eqnarray}\label{equ26}
 I&\geq &-\int_{\Omega_T}2(t+b) \sum_{i,k=1}^n a_{ik}\partial_{x_i} f_j \partial_{y_k} f_jdxdydt\nonumber\\
&&+\int_{\Omega_T}\sum_{i,k=1}^n \big(\frac12a_{ik}+(t+b)\partial_t a_{ik}\big)
\partial_{x_i} f_j \partial_{x_k} f_j dxdydt+\int_{\Omega_T}
 \frac{\alpha-\frac12}{2(t+b)}|f_j|^2dxdydt\nonumber \\
 &\geq& \int_{\Omega_T}\frac{1}{4\lambda}|\nabla_x f_j|^2-\lambda (t+b)|\nabla_y f_j|^2+
 \frac{\alpha-1}{2(t+b)}|f_j|^2.\quad
\end{eqnarray}

Note that \begin{eqnarray}\label{equ27} \int_{\Bbb R^m} |\partial_y
f_j|^2=\int_{\Bbb R^m} \eta^2|\psi(\eta/{2^j})\hat{f}(\eta)|^2\leq 4^{j+1}
\int_{\Bbb R^m} |\psi(\eta/{2^j})\hat{f}(\eta)|^2.\end{eqnarray}

Let $\alpha_0>1$, and we may assume that the parameter $\alpha$ in the function
$\phi(t)=(t+b)^{-\alpha}$ satisfies
$$
\alpha > \alpha_0>0.
$$
Moreover, we choose
$$j_0=max\{j\in {\Bbb Z};4^{j+1}\leq \frac{\alpha_0-1}{4(t+b)}, t\in (0,t_2]\}.$$
Hence when $\alpha>\alpha_0$ and $t_2\leq (16\lambda^2)^{-1}$, for
$j\leq j_0$, from (\ref{equ26}) and (\ref{equ27}), we obtain the following
Carleman-type inequality
\begin{eqnarray}\label{equ28}
\int_{\Omega_T}(t+b)|\phi L \Delta_j g|^2& \geq &\int_{\Omega_T}\frac{1}{4\lambda}|\nabla_x f_j|^2+
\frac{\alpha-1}{4(t+b)}|f_j|^2 \\
& =& \int_{\Omega_T}\frac{1}{4\lambda}(t+b)^{-2\alpha}|\nabla_x
\Delta_j g|^2+ \frac{\alpha-1}{4}(t+b)^{-2\alpha-1}|\Delta_j
g|^2.\nonumber
\end{eqnarray}
Then we finished the proof of Lemma 2.1. $\Box$

{\bf Proof of Theorem \ref{result01}:}
Since $g=\chi(t)u$, then
$$L\Delta_j g=\Delta_j (\chi'u+\chi c(x,t)u+\chi d(x,t)\cdot \nabla_x u).$$
By the above Carleman inequality (\ref{equ28}), we deduce
\begin{eqnarray}\label{equ29}
&\int_{\Omega_T}(t+b)^{-2\alpha+1}|(\chi'\Delta_j u+\chi c(x,t)\Delta_j u+
\chi d(x,t)\cdot \nabla_x\Delta_j u)|^2\nonumber
\\ \geq& \int_{\Omega_T}\frac{1}{4\lambda}(t+b)^{-2\alpha}|\chi  \nabla_x\Delta_j u|^2+
\frac{\alpha-1}{4}(t+b)^{-2\alpha-1}|\chi \Delta_j u|^2,
\end{eqnarray}
where we used the assumption $|c(x,t)|, |d(x,t)|\leq \lambda$
and the choice of $t_2$ which satisfies $t_2\leq (16\lambda^2)^{-1}
$. Consequently
$$
\int_{\Omega_T}(t+b)^{-2\alpha+1}|(\chi'\Delta_j u)|^2\geq
\int_{\Omega_T}\frac{\alpha-1}{8}(t+b)^{-2\alpha-1}|\chi \Delta_j
u|^2.
$$
Summing for all $j\leq j_0$, by the inequality ($\ref{equ23}$), we
obtain
$$
\int_{\Omega_T}(t+b)^{-2\alpha+1}|(\chi' u)|^2\geq
C(K)\int_{\Omega_T}\frac{\alpha-1}{8}(t+b)^{-2\alpha-1}|\chi
S_{j_0-1} u|^2.\eqno(2.8)
$$
Then
$$
\int_0^{t_1}\int_{R^{n+m}}|S_{j_0-1} u|^2\leq
\frac{16}{C(K)(\alpha-1)(t_2-t_1)^2}\int_{t_1}^{t_2}\int_{R^{n+m}}
\frac{(t_1+b)^{2\alpha+1}}{(t+b)^{2\alpha-1}}|u|^2. \eqno(2.9)$$

Now we let $\alpha\rightarrow\infty$ in (2.9), then we obtain
$$S_{j_0-1} u\equiv0.$$ And then $u\equiv0$ in $R^2\times (0,t_1)$ by
the choice of $j_0$ and $\alpha_0\rightarrow\infty$. Again we obtain
that $u\equiv0$ in $R^2\times (0,t_2)$ since $t_1$ can approach
$t_2$. Finally, we have $u\equiv 0$ in $R^2\times (0,T)$ after the
iteration, since $t_2$ only depends on L. Hence we have completed the proof
of Theorem 1.1. $\Box$

\begin{rmk}
The assumption that the coefficients $a_{ik}$, $c$ and $d$ are
independent of $y$ seems to be only a technique assumption. However,
we do not know how to remove it in general. The main difficulty for
the Kolmogorov operator in our case is the lose of derivative
estimates in $y$ direction. On the other hand, the recent regularity
result (see \cite{WZ}) shows that one can recover the regularity even in
$y$ direction.
\end{rmk}

Here we give another proof by using the frequency function method
(see [PO], [EKPV]). We consider the differential inequality
$$
|(\sum_{i,k=1}^n\partial_{x_i}(a_{ik}(x,t)\partial_{x_k})u
+\sum_{l=1}^mx_l\partial_{y_l}u+\partial_t u|\leq
\lambda(|u|+|\nabla_x u|).\eqno(2.10)
$$

 Let
$$e(t)=\int_{{\Bbb R}^{n+m}}u^2dx dy,$$
$$d(t)=\int_{{\Bbb R}^{n+m}}(\sum_{i,k=1}^na_{ik}(x,t)\partial_{x_i}u\partial_{x_k}u)dx dy,$$
and $$h(t)=\frac{d(t)}{e(t)}.$$ Then we have the following monotonicity inequality lemma.

\begin{lem}
Suppose that $u$ satisfies (2.10). In addition to the condition of
Theorem 1.1, we assume that for some constant $M$
$$
\int_{{\Bbb R}^{n+m}}|\nabla_yu|^2\leq M\int_{{\Bbb R}^{n+m}}|u|^2.\eqno(2.11)
$$
Then there exits a constant $C=C(\lambda, M)$ such that
$$\dot{h}(t)\geq -C(\lambda,
M)[h(t)+1].\eqno(2.12)$$
\end{lem}
\begin{rmk}
This is the corresponding monotonicity inequality for the Kolmogorov
operator. There is an additional assumption (2.11) which seems
necessary in our approach. And in application, we again need to make
use of the Littlewood-Paley decomposition in $y$ direction.
\end{rmk}

{\it Proof:} By our assumption, one can calculate directly
\begin{eqnarray*}\dot{e}(t)=2\int_{{\Bbb R}^{n+m}} uu_t
=\int_{{\Bbb R}^{n+m}} u Lu+2\int_{{\Bbb R}^{n+m}} u(\partial_t u+\sum_{l=1}^mx_l\partial_{y_l} u-\frac12
Lu),\end{eqnarray*}
and
$$d(t)=\int_{{\Bbb R}^{n+m}} u(\partial_t u+\sum_{l=1}^mx_l\partial_{y_l} u-\frac12 Lu)-\frac12 \int_{{\Bbb R}^{n+m}} u Lu.$$
Hence
$$
\dot{e}(t)d(t)=2[\int_{{\Bbb R}^{n+m}} u(\partial_t u+\sum_{l=1}^mx_l\partial_{y_l} u-\frac12
Lu)]^2-\frac12 [\int_{{\Bbb R}^{n+m}} u Lu]^2.\eqno(2.13)
$$
By our assumption on the coefficient $a_{ik}(x,t)$, we have
\begin{eqnarray*}&&\dot{d}(t)\\&=&\int_{{\Bbb R}^{n+m}} \sum_{i,k=1}^n\partial_ta_{ik}(x,t)\partial_{x_i}u\partial_{x_k}u+2\int_{{\Bbb R}^{n+m}} \sum_{i,k=1}^na_{ik}(x,t)\partial_{x_i}u\partial_{x_k}\partial_tu\\
&=&\int_{{\Bbb R}^{n+m}} \sum_{i,k=1}^n\partial_ta_{ik}(x,t)\partial_{x_i}u\partial_{x_k}u-2\int_{{\Bbb R}^{n+m}} \sum_{i,k=1}^n\partial_{x_k}(a_{ik}(x,t)\partial_{x_i}u)\partial_tu\\
&=&\int_{{\Bbb R}^{n+m}} \sum_{i,k=1}^n\partial_ta_{ik}(x,t)\partial_{x_i}u\partial_{x_k}u+2\int_{{\Bbb R}^{n+m}} (\partial_t u+\sum_{l=1}^mx_l\partial_{y_l} u)(\partial_t u+
\sum_{l=1}^mx_l\partial_{y_l}- Lu)dx dy\\
&&-2\int_{{\Bbb R}^{n+m}} \sum_{i,k=1}^na_{ik}(x,t)\partial_{x_i}u \partial_{y_k}u\\
&=&\int_{{\Bbb R}^{n+m}} \sum_{i,k=1}^n\partial_ta_{ik}(x,t)\partial_{x_i}u\partial_{x_k}u+2\int_{{\Bbb R}^{n+m}} (\partial_t u+\sum_{l=1}^mx_l\partial_{y_l}  u-\frac12 Lu)^2\\
&&-\frac12\int_{{\Bbb R}^{n+m}}|Lu|^2-2\int_{{\Bbb R}^{n+m}} \sum_{i,k=1}^na_{ik}(x,t)\partial_{x_i}u \partial_{y_k}u\\
&\geq&2\int_{{\Bbb R}^{n+m}}(\partial_t u+\sum_{l=1}^mx_l\partial_{y_l} u-\frac12
Lu)^2-\frac12\int_{{\Bbb R}^{n+m}}|Lu|^2 -C(\lambda,M)[d(t)+e(t)].
\end{eqnarray*}
Hence together with (2.13) we obtain
$$e(t)\dot{d}(t)-d(t)\dot{e}(t)\geq -C(\lambda,M)[d(t)+e(t)]-\frac12\int_{{\Bbb R}^{n+m}}|Lu|^2.\eqno(2.14)$$
Consequently, we have the following monotonicity inequality
$$\dot{h}(t)=\frac{\dot{d}(t)e(t)-\dot{e}(t)d(t)}{e(t)^2}\geq -C(\lambda, M)
[h(t)+1]-\frac{\int_{{\Bbb R}^{n+m}}|Lu|^2}{2 e(t)}.$$ Together with (2.10),
then (2.12) follows easily. Then we proved Lemma 2.2. $\Box$

{\bf Proof of  Theorem 1.1 by the monotonicity inequality:} \\
We assume that $e(t)\equiv 0$ for $0\leq t
\leq t_1$, and $e(t)>0$ as $t_1<t\leq T$. For $t_1<t<t_2\leq T$, we
integrate the above inequality from t to $t_2$
$$\log{\frac{h(t_2)+1}{h(t)+1}}\geq -C(\lambda, M)T,$$
which yields $$ h(t)\leq C(\lambda,M, T, h(t_2)).\eqno(2.15)$$ Since
$$\frac{\dot{e}(t)}{e(t)}\leq 3h(t)+C(\lambda)\leq C(\lambda,M, T,
h(t_2)),$$ integrating from t to $t_2$, we have
$$e(t_2)\leq e(t)C(\lambda,M, T, h(t_2)).\eqno(2.16)$$
Let $t\rightarrow t_1$, we get $e(t_2)=0$ which is a contradiction.
Using the same Littlewood-Paley decomposition as Theorem 1.1, we
could replace u with $\triangle_j u$, since $\triangle_j u$
satisfies the inequality (2.11). Then the remaining arguments are
similar to that of the previous proof of Theorem 1.1. We complete the proof
of our theorem by the frequency function method. $\Box$

\end{document}